\documentclass[preprint,12pt]{elsarticle}

\usepackage{hyperref,empheq}
\usepackage[active]{srcltx}

\usepackage[T1]{fontenc}
\usepackage[utf8]{inputenc}
\usepackage[a4paper]{geometry}
\usepackage{hyperref,empheq}
\usepackage[active]{srcltx}
\usepackage{amsmath,amssymb}
\newtheorem{theorem}{Theorem}
\newtheorem{proposition}{Proposition}

\newtheorem{definition}{Definition}
\newtheorem{lemma}{Lemma}

\newcommand{\R}{\mathbb{R}}

\renewcommand{\L}{\mathbb{L}}
\newcommand{\noi}{\noindent}
\renewcommand{\phi}{\varphi}
\renewcommand{\d}{\mathrm{d}}
\newcommand{\veps}{\varepsilon}

\journal{Jounal of Geometry and Physics}

\begin{document}

\begin{frontmatter}


 \tnotetext[l5]{Work supported in part by FMJH Program PGMO and from the support of
 	EDF-Thal\`es-Orange (PGMO grant no. 2016-1753H)}

\title{Non-integrability\\
		of the minimum-time Kepler problem\tnoteref{5}}


\author[1]{M.~Orieux}

\author[2]{J.-B.~Caillau}

\author[3]{T.~Combot}

\author[4]{J.~F\'ejoz}

\address[1]{CEREMADE, Univ.\ Paris Dauphine,
	Place du Mar\'echal de Lattre de Tassigny, F-75016 Paris
	(\texttt{orieux@ceremade.dauphine.fr})}

\address[2]{LJAD, Univ.\ C\^ote d’Azur \& CNRS/Inria, Parc Valrose, F-06108 Nice
	(\texttt{caillau@unice.fr})}

\address[3]{Institut math., UBFC \& CNRS, 9 avenue Savary, F-21078 Dijon
	(\texttt{thierry.combot@ubfc.fr})}

\address[4]{CEREMADE, Univ.\ Paris Dauphine,
	Place du Mar\'echal de Lattre de Tassigny, F-75016 Paris \& IMCCE, Observatoire de Paris, 77 avenue Denfert Rochereau, F-75014 Paris
	(\texttt{jacques.fejoz@dauphine.fr})}

\begin{abstract}
 We prove, using Moral\`es-Ramis theorem, that the minimum-time controlled Kepler problem is not
 meromorphically integrable in the Liouville sens on the Riemann surface of its
 Hamiltonian.

\end{abstract}


\end{frontmatter}

\section{Introduction}
\noi The Kepler problem
\begin{equation}
\label{eq:Kepler}
\ddot{q}+\frac{q}{\|q\|^3}=0, \quad q \in \R^2\setminus\{0\}.
\end{equation}
is a classical reduction of the two-body
problem~\cite{arnol2013mathematical}. Here, we think of $q$ as the position of a 
spacecraft, and of the attraction as the action of the Earth. We are
interested in controlling the transfer of the spacecraft from one
Keplerian orbit towards another, in the plane. Denoting
$v = \dot q$ the velocity, and the adjoint variables of $q$ and $v$ by $p_q$ and
$p_v$, the minimum time dynamics is a Hamiltonian system with
\begin{equation}
\label{eq:0}
H(q,v,p_q,p_v)=p_q.v-\frac{p_v.q}{\|q\|^3}+\|p_v\|,
\end{equation}
as is explained in section~\ref{sec:reminder}.
Prior studies of this problem can be found in \cite{caillau2012minimum, caillau2001}.
The controlled Kepler problem can be embedded in the two parameter
family obtained when considering the control of the circular restricted three-body
problem:
\begin{equation} \label{eq51}
\ddot{q}+\nabla_q\Omega_\mu(t,q) = \veps u,
\end{equation}
where
\begin{eqnarray*}
	\Omega_\mu(t,q) &=&
	-\frac{1-\mu}{\sqrt{(q_1+\mu\cos t)^2+(q_2+\mu\sin t)^2}}\\
	&& -\frac{\mu}{\sqrt{(q_1-(1-\mu)\cos t)^2+(q_2-(1-\mu)\sin t)^2}}
\end{eqnarray*}
is the potential parameterized by the ratio of masses, $\mu \in [0,1/2]$,
and where $u \in \R^2$ is the control, whose amplitude is modulated by the second
parameter, $\veps \geq 0$. Alternatively to time minimization, minimization of the
$\L^2$ norm of the control can be considered,
\[ \int_0^{t_f} u^2(t)\,\d t \to \min. \]
This is the so-called energy cost. In the uncontrolled model ($\veps=0$), it is
well known that the Kepler case ($\mu=0$) is integrable and geodesic (there exists a Riemannian metric such that Keplerian curves are geodesics of this metric \cite{moser-1970a,osipov-1977a})
while there are obstructions to integrability for positive $\mu$.
In the controlled case ($\veps>0$), the Kepler
problem for the energy cost has been shown to be integrable (and geodesic)
when suitably averaged (see \cite{caillau2009} for a survey). The aim of this paper is
to study the integrability properties of the Kepler problem for time minimization.

The pioneering work of Ziglin in the 80's \cite{ziglin1982branching},
followed by the modern formulation of differential Galois theory in
the late 90's by Moral{\`e}s, Ramis and Sim{\'o}
\cite{morales2001galoisian,morales2007integrability}, have led to a
very diverse literature on the integrability of Hamiltonian systems.
According to Pontrjagin's Maximum principle, one can turn general
optimization problems with dynamical constraints into Hamiltonian
systems, which are generally not everywhere differentiable. Optimal
control theory thus provides an abundant class of dynamical systems
for which integrability is a central question. Yet,
differential Galois theory has not so often been applied in this context
(see, \emph{e.g.}, \cite{combot-2015a}),
in part because of the difficulty brought by the singularities.
Notwithstanding theses singularities (vanishing of the adjoint variable $p_v$, here),
we show how to apply these ideas to the system~\eqref{eq:0}.

\section{Setting}

\subsection{The minimum time controlled Kepler problem}
\label{sec:reminder}

We first recall some classical facts on optimal control. We refer for
example to the book of Agrachev and Sachkov \cite{agrachev13} for more details.
Let $M$ be an $n$-dimensional smooth manifold and $U$ an arbitrary
subset of $\mathbb{R}^m$ (typically a submanifold with boundary).
A controlled dynamical system is a smooth family of vector fields
\[f : M\times U\rightarrow TM\]
parameterized by the control values.
Admissible controls are measurable functions valued in the subset $U$.
A preliminary question is the following:
Is some final state $x_f$ accessible from some initial state $x_0$,
\emph{i.e.} does the system
\[ \dot{x}(t)=f(x(t),u(t)),\quad u(t) \in U, \]
\[ x(0)\ =x_0,\quad x(t_f)=x_f, \]
have a solution for some admissible control?  The system
is said to be controllable if the answer is positive for all
possible initial and final states $x_0,x_f\in M$.
The controlled Kepler problem, associated with~\eqref{eq:Kepler}, is
\[ \ddot{q}+\frac{q}{\|q\|^3}=u, \quad q \in \R^2\setminus\{0\}, \quad
u_1^2+u_2^2\leq1, \]
\[ (q(0),\dot q(0))\ \ =(q_0,v_0),\quad
(q(t_f),\dot q(t_f))=(q_f,v_f), \]
where $q$ is the position vector of a spacecraft and where
the control $u$ is the thrust of the engine. The
thrust is obviously bounded; here we assume that it is valued in the
Euclidean unit ball. (Note that, with respect to (\ref{eq51}),
we have chosen $\veps=1$; as will be clear from Section~\ref{sproof},
this does not restrict the generality of the analysis.) 

\begin{proposition}[\cite{caillau2001}]
	The Kepler problem is controllable.
\end{proposition}

\noi This is a consequence of two facts: The Lie algebra generated by the
drift and the vector field supporting the control generate the whole
tangent space at each point (which entails some local
controllability), and the uncontrolled flow (or \emph{drift}) of the
Kepler problem is recurrent. Under some additional convexity and compactness assumptions, one is
then able to retrieve existence of optimal controls.

We now deal with such optimal controls.
We restrict ourselves to integral cost functions, that is to problems of the form
%
\begin{equation}
\label{eq:opti}
\begin{cases}
\dot{x}(t)=f(x(t),u(t)),\\
x(0)\ =x_0,\quad x(t_f)=x_f,\\
\int_{0}^{t_f} L(x(t),u(t))\, \text{dt}\rightarrow \min
\end{cases}
\end{equation}
%
where the final time $t_f$ can be fixed or not, and $L:M\times U\rightarrow\mathbb{R}$ 
is a smooth function.
In the early 60's, Pontrjagin and his coauthors realized that necessary conditions
for optimality could be stated in Hamiltonian terms. By $T^*M$ we denote the cotangent
bundle of the manifold $M$.

\begin{definition} \label{def1}
	The associated pseudo-Hamiltonian is
	\[H: T^*M\times\R\times U \to \R, \quad (x,p,p^0,u) \mapsto \langle
	p,f(x,u)\rangle+p^0 L(x,u).\]
\end{definition} 

\noi The following fundamental result is Pontrjagin Maximum Principle
\cite{pontryagin1987mathematical}  (see
\cite{agrachev13} for a modern presentation).

\begin{theorem} \label{th1}
	If $(x,u)$ solves \eqref{eq:opti}, there exists a Lipschitzian function
	$p(t)\in T^*_{x(t)}M$, $t\in[0,t_f]$, a constant $p^0\leq 0$,
	$(p(t),p^0)\neq 0$, such that, almost everywhere,
	\begin{itemize}
		\item[(i)] $(x,p)$ is a solution of the Hamiltonian system associated with
		$H(\cdot, \cdot, u(t))$:
		\[ \dot{x}=\frac{\partial H}{\partial p}(x,p,u),\quad
		\dot{p}=-\frac{\partial H}{\partial x}(x,p,u), \]
		\item[(ii)] $H(x(t),p(t),u(t))=\max_{v \in U}H(x(t),p(t),v)$.
	\end{itemize}
	Such curves $(x,p)$ are called extremals. As a consequence of the
	maximization condition, the pseudo-Hamiltonian evaluated along an extremal is
	constant. Moreover, if the final time is free then this constant is zero.
\end{theorem}

\noi This powerful result has some downsides. The Hamiltonian is
defined on the cotangent bundle of the original phase space, and thus
the dimension is doubled. Besides, the maximization condition, which "eliminates the
control" and allows to obtain a truly Hamiltonian system in $(x,p)$ only, might generate
singularities (that is non-differentiability points of the maximized Hamiltonian which is in general only Lipschitzian as a function of time when evaluated along an extremal).
The above theorem applies to time minimization with $L\equiv 1$ (and free final
time). In this case, the non-positive constant $p^0$ is only related to the level of the
Hamiltonian, and we will not mention it in the sequel as we will not discuss the
implications of having \emph{normal} ($p^0 \neq 0$) or \emph{abnormal} ($p^0=0$) extremals.

\subsection{Main result}
The minimum time Kepler problem can be stated according to
\begin{equation}\label{eq:mt}
\begin{cases}
\ddot{q}+\frac{q}{\|q\|^3} =u,\quad \|u\|\leq 1,\\
(q(0),\dot q(0)) = (q_0,v_0),\quad (q(t_f),\dot q(t_f))  =  (q_f,v_f),\\
t_f\rightarrow \min,
\end{cases}
\end{equation}
where, as before, $q \in \R^2$ is the position vector and $u \in \R^2$
the control. It will be convenient to use the same notations as in the
general problem~\eqref{eq:opti} and let
\[q=(x_1,x_2), \quad \dot q=(x_3,x_4),\]
be the coordinates on the initial phase space
$M = (\R^2 \setminus \{0\}) \times \R^2 $. According to Definition~\ref{def1},
the pseudo-Hamiltonian is then
\begin{equation} \label{eq20}
H(x,p,u) = p_1x_3+p_2x_4 - \frac{p_3x_1+p_4x_2}{(x_1^2+x_2^2)^{3/2}}
+ p_3u_1+p_4u_2.
\end{equation}
According to Theorem~\ref{th1}, minimizing trajectories must be projections on $M$ of
integral curves of the Hamiltonian that has to be maximized over the unit disk.
The maximized Hamiltonian is readily equal to 
\[H(x,p) = p_1x_3+p_2x_4 - \frac{p_3x_1+p_4x_2}{(x_1^2+x_2^2)^{3/2}} +
\sqrt{p_3^2+p_4^2}\]
on $T^*M$, while the control is given by
\[ u = \frac{1}{\sqrt{p_3^2+p_4^2}}(p_3,p_4) \]
whenever $p_3$ and $p_4$ do not vanish simultaneously. 
Now, let 
\[\mathcal{M}=\{(x,p,r)\in\mathbb{C}^8\times \mathbb{C}_*^2,\;
r_1^2=x_1^2+x_2^2,\;\; r_2^2=p_3^2+p_4^2\}\]
be the Riemann surface of $H$. It is a complex symplectic manifold
(with local Darboux coordinates $(x,p)$ outside the singular
hypersurface $r_1r_2=0$), over which $H$ extends meromorphically, and
even rationally, since
\begin{equation}
\label{eq:H}
H(x,p,r)=p_1x_3+p_2x_4-\frac{p_3x_1+p_4x_2}{r_1^3}+r_2.
\end{equation}
The Hamiltonian $H$ has four degrees of freedom, hence (see
\cite{arnol2013mathematical}) the meromorphic Liouville integrability
of $H$ over $\mathcal{M}$ would mean that there would exist three independent
first integrals, in addition to $H$ itself, almost everywhere in
$\mathcal{M}$. The aim of this paper is to prove that it is not the case.

\begin{theorem}
	The minimum time Kepler problem is not meromorphically Liouville
	integrable on $\mathcal{M}$.
	\label{thm:ni}
\end{theorem}

\noi It is well known that the classical Kepler problem is integrable, and
even super integrable (since there are more first integrals than
degrees of freedom, as a result of Kepler's first law and of the
dynamical degeneracy of the Newtonian potential---see
for instance \cite{Fejoz:2004a}).
On the opposite, the three-body problem is
not as is known after the seminal work of Poincar{\'e} (for recent accounts on this topic see, \emph{e.g.}, \cite{Combot:2012:nbp, Julliard:2000:Bruns,
	Poincare:1892, Tsygvintsev:2003}). Similarly, the above theorem
asserts that lifting the Kepler problem to the cotangent bundle and
introducing the singular control term $r_2$ breaks integrability.

This result prevents the existence of enough complex analytic (and even
meromorphic) first integrals to ensure integrability over
$\mathcal{M}$. Or course, it does not prevent the existence of an
additional real first integral
which would have a natural frontier asymptotic to the 
real domain and thus, would not extend to the complex plane.
%
Future work might be dedicated to investigate either or not
Theorem \ref{thm:ni} holds for real first integrals.

\section{Proof of Theorem~\ref{thm:ni}} \label{sproof}
\noi The rest of the article is devoted to proving the theorem. Our proof
consists in studying the variational equation along some integral
curve of~\eqref{eq:H}. In order to carry out this computation, we
choose a collision orbit, with the drawback that it requires some
regularization. We also note that there exist effective tools to perform this kind of
computations (see, \emph{e.g.}, \cite{weil-2012a}).
The algebraic obstruction to Liouville integrability
comes from the theorem below of Moral{\`e}s and Ramis, which we now recall. We follow the presentation of Singer in \cite{singer2009}.

\subsection{Some facts of Galois differential theory}
Consider a linear differential equation $(L) :Y'=AY$, $A\in M_n(k)$,
$k$ being a differential field whose field of constants $k_0$ is algebraically closed, and of characteristic zero. We want the Galois group to be the group of symmetries preserving all algebraic and differential relations of a basis of solutions. We consider the  polynomial ring
$$S=k[Y_{1,1},\dots,Y_{n,n},1/\hbox{det}(Y)]$$
where $Y$ is an $n\times n$ matrix. This ring has a derivation provided by the differential system $Y'=AY$. We now consider a maximal differential ideal $M$ of $S$, and the quotient $R=S/M$. This quotient satisfies the following

\begin{definition}[Picard-Vessiot field]
	A Picard-Vessiot ring for $Y' = AY$ is a differential ring $R$ over $k$ such that
	\begin{itemize}
		\item[(i)] The only differential ideals of $R$ are $(0)$ and $R$.
		\item[(ii)] There exists a fundamental matrix $Z\in GL_n(R)$ for the equation
		$Y' = AY$.
		\item[(iii)] $R$ is generated as a ring by $k$, the entries of $Z$ and $1/\hbox{det}(Z)$ . 
	\end{itemize}
\end{definition}
It turns out that the choice of the maximal differential ideal $M$ always gives the same Picard-Vessiot ring up to isomorphism. This ring is also a domain, thus allowing to consider the quotient field, the Picard-Vessiot field.

\begin{definition}[Galois group]
	The differential Galois group of $R$ over $k$ is the group of differential automorphism of $R$ preserving $k$, noted $Gal(R/k)$.
\end{definition}

\noindent For a differential system $Y'=AY$, if there is no ambiguity on the base field $k$. (For the case treated in this paper, the base field $k$ is $\mathbb{C}(z)$.) 
Given a fundamental matrix of solution $Z$ and a Galois group element $\sigma$, we have $Z'=AZ$, and thus applying $\sigma$, we also have $\sigma(Z)'=A\sigma(Z)$. Thus $\sigma(Z)$ is also a matrix of solutions; there exists a constant matrix $C$ such that $\sigma(Z)=ZC$, and as $\sigma$ is an automorphism, $C$ has to be invertible. So $Gal(R/k)$ can be represented as a group of $n\times n$ matrices.

\begin{proposition}
	The Galois group $Gal(R/k)\subset GL_n(k_0)$ is a linear algebraic
	group, i.e. the zero set in $GL_n(k_0)$ of a system of polynomials
	over $k_0$ in $n^2$ variables.
\end{proposition}

\noindent\textbf{Proof.}
	This can be obtained by letting a Galois group element $\sigma$ act (right multiplication by a matrix) on the differential ideal $I=(f_1,\dots,f_p)$. We can moreover assume that $f_i\in k[Y]$. As this does not change the degrees in the $Y_{i,j}$ and since $I$ must be stabilized, $\sigma(f_i)$ must belong to $I\cap k_{\max(\deg f_1,\dots, \deg f_p)}[Y]$. This condition is a condition of membership to a vector space, which provides algebraic conditions on the entries of the matrix $\sigma$.

\begin{proposition}[Fundamental Theorem of Differential Galois
	Theory]
	Let $K$ be a Picard-Vessiot field with differential Galois group $G$ over $k$.
	\begin{itemize}
		\item[(i)] There is a one-to-one correspondence between Zariski-closed subgroups
		$H\subset G$ and differential subfields $F$, $k\subset F\subset K$, given by
		$$H \subset G \rightarrow K^H = \{a \in K, \sigma (a) = a\; \forall \sigma \in H\}$$
		$$F \rightarrow Gal(K/F)= \{\sigma \in G, \sigma(a)=a \; \forall a\in F \}$$
		\item[(ii)] A differential subfield $F$, $k\subset F \subset K$, is a Picard-Vessiot extension of $k$ if and only if $Gal(K/F)$ is a normal subgroup of $G$, in which case $Gal(F/k) \simeq G/Gal(K/F)$.
	\end{itemize}
\end{proposition}

We are interested in non integrability for nonlinear Hamiltonian systems. The link with the Hamiltonian world is given by the celebrated theorem of 
Moral\'es-Ramis below. We recall that an algebraic group $G$ is said to be virtually Abelian if 
its connected component containing the identity is an Abelian subgroup of $G$.

\begin{theorem}[Moral{\`e}s-Ramis \cite{morales2001galoisian}]
	\label{thm:mr}
	Let $H$ be an analytic Hamiltonian on a complex analytic
	symplectic manifold and $\Gamma$ be a non constant solution. If $H$ is
	integrable in the Liouville sense with meromorphic first integrals,
	then the first order variational equation along $\Gamma$ has a
	virtually Abelian Galois group over the base field of meromorphic
	functions on $\Gamma$.
\end{theorem}

\noi The main idea behind this theorem is that if $H$ is Liouville integrable, then so are the linearized equations near a non constant solution $\Gamma$. More precisely, thanks to Ziglin's Lemma below, the first integrals of $H$ can be transformed in such a way that their first non trivial term in their series expansion near $\Gamma$ are functionally independent.

\begin{lemma}[Ziglin's Lemma]
	Let $\Phi_1,\dots,\Phi_r\in k(x_1,\dots,x_n)$ be functionally independent functions. We consider $\Phi_1^0,\dots,\Phi_r^0$ the lowest degree homogeneous term for some fixed positive weight homogeneity in $x_1,\dots,x_n$. Assume $\Phi_1^0,\dots,\Phi_{r-1}^0$ are functionally independent. Then there exists a polynomial $\Psi$ such that the lowest degree homogeneous term $\Psi^0$ of $\Psi(\Phi_1,\dots,\Phi_r)$ is such that $\Phi_1^0,\dots,\Phi_{r-1}^0,\Psi^0$ are functionally independent.
\end{lemma}

\noindent Applying this Lemma recursively, we prove that if a Hamiltonian system admits a set of commuting, functionally independent meromorphic first integrals on a neighbourhood of a curve, then their first order terms, after possibly polynomial combinations of them, are also commuting, functionally independent meromorphic first integrals of the linearized system along the curve.
Moral{\`e}s-Ramis \cite{morales2001galoisian} precisely proved that symplectic linear differential systems having such first integrals have a Galois group whose identity component is Abelian. This result can be expected knowing that the Galois group leaves invariant every first integral, so the more first integrals, the smaller the Galois group.

We will need the definition of the monodromy group of a
linear differential equation. Let us consider a differential system $Y'=AY,\; A\in M_n(\mathbb{C}(x))$. We note $S=\mathbb{P}^1 \setminus \{\hbox{singularities of } A\} $.
Let us consider a point $z_0\in S$ and a closed oriented curve $\gamma\subset S$, with $x_0 \in\gamma$. There exists a basis of solutions $Z$ on a neigbourhood of $x_0$, holomorphic in $z$. We now use analytic continuation along the loop $\gamma$ to extend this basis of solutions. However, it cannot \emph{a priori} be extended to a whole neighbourhood of $\gamma$, because after one loop, the basis of solutions $Z_\gamma$ at $x_0$ could be different. This defines a matrix $D_\gamma \in GL_n(\mathbb{C})$ such that $Z_\gamma =Z D_\gamma$ and thus a homomorphism
$$\hbox{Mon} :\pi_1(S,x_0) \rightarrow GL_n(\mathbb{C}), \quad \hbox{Mon}(\gamma)= D_\gamma.$$
This homomorphism carries the group structure of $\pi_1(S,x_0)$, and thus its image is also a group.

\begin{definition}
	The image of the application $\hbox{Mon}$ is called the monodromy group.
\end{definition}

\noindent Note that the monodromy group depends on the choice of $Z$, so it is only
determined up to conjugation. Since analytic continuation preserves
analytic relations, the monodromy group is a subset of
the differential Galois group over the base field of meromorphic functions on $S$; in particular, it is included in the differential Galois group over the base field of rational functions. For Fuchsian systems (all singularities are regular singularities, i.e. the growth at singularities of solutions is at most polynomials), we have moreover the following.

\begin{theorem}[Schlesinger density theorem \cite{sauloy16}]\label{thmschles}
	Let $(E):\ Y'=AY$ be a Fuchsian differential linear equation with coefficients in $\mathbb{C}(x)$ and let
	$\Pi$ be its monodromy group. Then $\Pi$ is dense for the Zariski
	topology in the Galois group of the Picard-Vessiot extension of
	$(E)$ over the base field of rational functions: $\overline{\Pi} = \textrm{Gal}(A)$.
\end{theorem}

\subsection{A collision orbit}

In order to find an explicit solution of \ref{eq20}, let us define the
$4$-dimensional symplectic submanifold
\[ S = \{(x,p,r)\in \mathcal{M}\ |\ x_2=x_4=p_2=p_4=0,r_1=x_1,r_2=-p_3\}. \]
As $S$ is the phase space of the controlled Kepler problem on the line
(collision orbit)
parameterized by $q_1$, it is invariant. On the interior of $S$,
$(x_1,x_3,p_1,p_3)$ is a set of (Darboux) coordinates and, in
restriction to $S$, the Hamiltonian reduces to
\[H(x,p)=p_1x_3-\frac{p_3}{x_1^2}-p_3,\]
so the Hamiltonian vector field on $S$ is
\[
\begin{cases}
\dot x_1 = x_3\\
\dot x_3 =- 1 -\frac{1}{x_1^2}\\
\dot p_1 = - \frac{2p_3}{x_1^3}\\
\dot p_3 = - p_1.
\end{cases}
\]
In particular,
\begin{equation}
\label{eq:x1p3}%
\begin{cases}
\ddot{x}_1=-1-\frac{1}{{x_1}^2} \\
\ddot{p}_3-\frac{2p_3}{x_1^3}=0.
\end{cases}
\end{equation}
As is known since the work of Charlier and Saint Germain on the Kepler problem with a
constant force (see \cite{beletsky01}), the function
\[C=\frac{1}{2}x_3^2+x_1-\frac{1}{x_1}\]
is a first integral on $S$ and $H_{|S}$ is integrable. Let us change time to $s=x_1(t)$ and denote by $'=\frac{d}{ds}$ the derivation with respect to this new time. It suffices to find an obstruction in this modified time, as
explained at the end of the proof.

Using~\eqref{eq:x1p3}, we see that the variable $p_3$ satisfies the linear differential equation
\[2\left(C+\frac{1}{x_1}-x_1\right)p_3''(x_1) -
\left(1+\frac{1}{{x}_1^2}\right)p_3'(x_1)-\frac{2p_3(x_1)}{x_1^3}=0,\] 
which yields
\[p_3(x_1) = \frac{\sqrt{-Cx_1+x_1^2-1}}{\sqrt{x_1}}
\left(c_1\int\frac{x_1^{3/2}}{(-Cx_1+x_1^2-1)^{3/2}}dx_1+c_2\right)\]
for some constants of integration $c_1$ and $c_2$. Here the symbol
$\int f(x_1)dx_1$ denotes some primitive of $f$ with respect to the
variable $x_1$. It suffices to find one particular integral curve along which the
variational equation has a non virtually Abelian Galois group. To this end, we 
consider the simple---but rich enough---case $c_1=0$, $c_2=1$.
\[p_3(x_1)=\frac{\sqrt{-Cx_1+x_1^2-1}}{\sqrt{x_1}}\cdot \]
Using the expression of the first integral $C$ and of the vector
field, we deduce
\[x_3(x_1)=\sqrt{2}\frac{\sqrt{-Cx_1+x_1^2-1}}{\sqrt{x_1}}, \;\;
p_1(x_1)=-\frac{1}{\sqrt{2}}\frac{x_1^2+1}{x_1^2}\cdot\]
Choosing $C=2i$ and some determination of the squares yields
a particularly simple solution $\Gamma$ drawn on
$S\subset \mathcal{M}$,
\begin{equation}\label{eqGamma}
\begin{cases}
x_1 = x_1,\\
x_2 = 0,\\
x_3 = \sqrt{2}\frac{x_1-i}{\sqrt{x_1}},\\
x_4 = 0,\\
\end{cases}
\quad 
\begin{cases}
p_1 = -\frac{x_1^2+1}{\sqrt{2} x_1^2},\\
p_2 = 0,\\
p_3 = \frac{x_1-i}{\sqrt{x_1}},\\
p_4 = 0.
\end{cases}
\end{equation}

\subsection{Normal variational equation}
In the initial time, the linearized equation along $\Gamma$ is the
Hamiltonian vector field associated with the Hamiltonian $DH$ along
$\Gamma$:
\[\dot{Z}(t)=A(t)Z(t), \quad A(t)=J \, D^2 H(\Gamma(t)),\]
where $J$ is the Poisson structure. In the coordinates
$(x_1,...,x_4,p_1,...,p_4)$, 
\[
J=\begin{pmatrix}
0_4 & I_4\\
-I_4 & 0_4
\end{pmatrix}.\]
We will keep on using time $x_1$, instead of the initial time $t$,
writing
\[Z'(x_1(t))=\frac{1}{x_3(t)}A(x_1(t))Z(x_1(t)).\]
Let us now reorder coordinates according to $(x_1,x_3,p_1,p_3,x_2,x_4,p_2,p_4)$.
Since $S$ is an invariant submanifold, the $8\times8$ matrix $A$ has
an upper triangular bloc structure
\[A=\begin{pmatrix}
A_1 & A_2\\
0 & A_3
\end{pmatrix}\]
with
\[A_3=\begin{pmatrix}
0 & 0 & 0 & \frac{1}{\sqrt{2}p_3}\\
-\frac{1}{\sqrt{2}p_3^2} & 0 & -\frac{1}{\sqrt{2}x_1^3p_3} & 0\\
0 & \frac{1}{\sqrt{2}p_3} & 0 & 0\\
-\frac{1}{\sqrt{2}x_1^3p_3} & 0 & \frac{3}{\sqrt{2}x_1^4} & 0
\end{pmatrix}.\]
Moral\`es-Ramis Theorem gives necessary conditions for Liouville
integrability in terms of the Galois group of this linear differential
system over the base field of meromorphic functions \emph{on}
$\Gamma$. Looking at the expression \eqref{eqGamma} of $\Gamma$, we
see that meromorphic functions on $\Gamma$ are just meromorphic
functions in $\sqrt{x_1}\in \mathbb{C}\setminus \{0,\pm \sqrt{i}\}$.
The block $A_3$ corresponds to infinitesimal variations in the normal
direction to $S$, which is the part where interesting phenomena might
occur.  As the Picard-Vessiot field is generated by all the components
of the solutions, the Picard-Vessiot field $K$ generated by the normal
variational equation
\[(L): \quad X'=A_3X, \quad X=(X_1,X_2,X_3,X_4)\]
is a subfield of the Picard-Vessiot field of the whole variational
equation, and thus $\hbox{Gal}(A) \supset \hbox{Gal}(A_3)$. That
$\hbox{Gal}(A_3)$ is not virtually Abelian will thus imply that
$\hbox{Gal}(A)$ itself is not virtually Abelian.
In order to reduce the system to a one dimensional linear equation, we
use the cyclic vector method on $A_3$: From $(L)$ we get
$X_1'=L_1(X_1,X_2,X_3,X_4)$, where $L_1$ is a linear form on $\R^4$,
thus by derivation,
\begin{align*}
	X_1''=&L_1(X_1',X_2,X_3,X_4) + L_1(X_1,X_2',X_3,X_4)\\
	+&L_1(X_1,X_2,X_3',X_4) + L_1(X_1,X_2,X_3,X_4')\\
	=&L_2(X_1,X_2,X_3,X_4).
\end{align*}
Iterating, we obtain
\begin{empheq}[left = \empheqlbrace]{align*}
X_1&=X_1,\\
X_1'&=L_1(X_1,X_2,X_3,X_4),\\
X_1''&=L_2(X_1,X_2,X_3,X_4),\\
X_1^{(3)}&=L_3(X_1,X_2,X_3,X_4),\\
X_1^{(4)}&=L_4(X_1,X_2,X_3,X_4).
\end{empheq}
The $L_i$'s are five linear forms on $\R^4$, so $X_1$ must satisfy some linear differential equation of order $4$ that we compute to be
\begin{multline} 
	X_1^{(4)}+
	\frac{2(3i-5x_1)}{x_1(i-x_1)} X_1^{(3)} +
	\frac{(-3x_1+i)(-29x_1+23i)}{4(x_1-i)^2x_1^2}X_1''\\
	-\frac{(i-3x_1)(7x_1+i)}{4(x_1-i)^2x_1^3}X_1'
	+\frac{3x_1+i}{4(x_1-1)^3x_1^4}X_1=0.
	\label{hyp}
\end{multline}
We find a solution of this equation of the form \small
\[X_1(x_1)=\frac{i-x_1}{\sqrt{x_1}}\left( c_1+c_2\int
\sqrt{x_1}(1+ix_1)^{-\frac{3}{2}-i\frac{\sqrt{3}}{2}}.{} _2{F}_1
(\gamma(x_1))dx_1 \right),\]
where ${}_2F_1$ is the Gauss hypergeometric function and
\[\gamma(x_1) =
\left(\frac{5}{2} - i\frac{\sqrt{3}}{2}, \frac{1}{2} +
i\frac{\sqrt{3}}{2}, 1+i\sqrt{3},1+i x_1\right).\]
The Picard-Vessiot field $K$ contains this solution and, as it is a
differential field, it also contains
\[\sqrt{x_1}(1+ix_1)^{-\frac{3}{2}-i\frac{\sqrt{3}}{2}} {}_2{F}_1(\gamma(x_1)).\] 
Noting $\tilde{K}$ the differential field generated by this function, we have $\tilde{K} \subset K$. Now the Galois group of ${}_2{F}_1(\gamma(x_1))$ over
$\mathbb{C}(x_1)$ is $SL_2(\mathbb{C})$ (see Kimura's table, \cite{kimura1969riemann}). By Galois correspondence, the Galois group of (\ref{hyp}) over the rational functions in $x_1$ admits $SL_2(\mathbb{C})$ as a subgroup.
The hypergeometric equation \eqref{hyp} is Fuchsian 
(all its singular points are regular), so thanks to Theorem~\ref{thmschles},
we know that its Galois group over the field of rational functions is the closure of its monodromy group.
Besides, the Galois group over meromorphic functions contains the
monodromy group, and of course, is included in the Galois group over
rational functions. Eventually, the Galois group of (\ref{hyp}) over meromorphic functions in
$x_1$ also contains $SL_2(\mathbb{C})$. Thus,
adding the algebraic extension $\sqrt{x_1}$, the Galois group can be
reduced to at most one subgroup of index $2$: The only possibility is that the identity component contains
$SL_2(\mathbb{C})$ again.
So the Galois group of $K$ over the base field of meromorphic functions
in $\sqrt{x_1}\in \mathbb{C}\setminus \{0,\pm \sqrt{i}\}$ contains
$SL_2(\mathbb{C})$ and is not virtually Abelian. According to Morales-Ramis,
this concludes the proof.

\end{document}